\documentclass[11pt]{article}

\usepackage{amssymb}
\usepackage{amscd}
\usepackage{cite}
\usepackage{epsfig}
\usepackage{rotate}

\newcommand{\beq}{\begin{equation}}
\newcommand{\eeq}{\end{equation}}
\newcommand{\Leq}[1]{\label{#1}\end{equation}}
\newcommand{\bdm}{\begin{displaymath}}
\newcommand{\edm}{\end{displaymath}}

\newtheorem{theorem}{Theorem}

\newtheorem{proposition}[theorem]{Proposition}

\newtheorem{conjecture}{Conjecture}

\def\B{{\cal B}}
\def\I{{\cal I}}
\def\M{{\cal M}}
\def\V{{\cal V}}
\def\F{\mathbb F}
\def\Q{\mathbb Q}
\def\Z{\mathbb Z}
\def\mod#1{\,({\rm mod\ }#1) }
\def\proof{\medskip\noindent {\bf Proof.}\hskip 10pt}
\def\endproof{\hfill\vbox{\hrule
        \hbox{\vrule\kern4pt\vbox{\kern4pt
        \kern4pt}\kern4pt\vrule}\hrule}\bigskip}

\title{Scaling dynamics of a cubic interval-exchange transformation}
\author{J. H. Lowenstein${}^\dag$ and F. Vivaldi}
\date{\it\small
${}^\dag$Dept.~of Physics, New York University, 2 Washington Place, New
York, NY 10003, USA\\
School of Mathematical Sciences, Queen Mary, University of London,
London E1 4NS, UK
}

\begin{document}
\renewcommand{\labelenumi}{$(\roman{enumi})$}
\maketitle
\begin{abstract}
We study the dynamics of renormalization of a specific interval exchange 
transformation which features exact scaling (the cubic Arnoux-Yoccoz model).
Using a symbolic space that describes both dynamics and scaling, we 
characterize the periodic points of the scaling map in terms of generalized
decimal expansions, where the base is the reciprocal of a Pisot number 
and the digits are algebraic integers.
The set of periodic points has a rich arithmetic and geometric 
structure: we establish rigorously some basic facts, and use 
extensive numerical experimentation to formulate a conjecture.
\end{abstract}
\vspace*{30pt}\centerline{\today}
\section{Introduction}\label{section:Introduction}

There are symbolic representations of dynamical systems 
which lead to interesting arithmetical phenomena. 
In such {\it arithmetic codings,} the symbols represent generalized 
`digits' in a number system which naturally describes the dynamics. 
For a survey of this rich area of research from a dynamical perspective, 
see \cite{Sidorov}; an arithmetical viewpoint is developed in 
\cite{Akiyama,AkiyamaBorbeliBrunottePethoThuswaldner,%
AkiyamaBrunottePethoThuswaldner}. 

In arithmetic codings ---even in the simplest settings--- one encounters 
at once interesting problems.
A well-known example is the doubling map of the unit interval (circle) 
$\gamma: \,x\mapsto 2x \mod{1}$, with the symbolic dynamics 
defined by the partition of the interval into two equal halves.
The symbols representing an orbit are just the binary digits of its 
initial point, so that the eventually periodic symbolic sequences 
correspond to the rational numbers, and the periodic sequences to the
rationals with odd denominator.
For each $n$, the periodic set 
\begin{equation}\label{eq:Lattice}
\mbox{Fix}(\gamma^n)=\frac{1}{2^n-1}\Z \cap [0,1)
\hskip 40pt n\geq 1
\end{equation}
is a lattice. The exponentially large denominator determines the 
maximal complexity (height) an $n$-cycle can have. 

The question of minimal complexity is altogether more difficult: 
when $2^n-1$ has a large number of divisors, cancellation becomes 
likely, leading to $n$-cycles with much smaller denominator. 
This phenomenon becomes evident when one looks at the twin asymptotic 
problem, namely that of determining the period of rational points 
with a given odd denominator. In this setting, as the denominator 
increases, the low-complexity cycles are those that are found first.
It is possible to show that this problem is closely related to the
so-called Artin's conjecture on primitive roots \cite{RamMurty}; 
at present, the asymptotic behaviour of the period of a rational point 
can only be established assuming the validity of the generalized Riemann 
hypothesis \cite{Moller}.

In the study of periodic orbits, a natural generalization of 
the above construction occurs in hyperbolic toral automorphisms. 
In place of $\Q$, we now have the field $\Q(\lambda)$, where
$\lambda$ is an eigenvalue of the automorphism, while the
lattice (\ref{eq:Lattice}), which is a $\Z$-module\footnote{\
given a ring $R$, an $R$-module is an an additive group equipped 
with multiplication by elements of $R$ ---see \cite{HartleyHawkes}}, 
is replaced by a $\Z[\lambda]$-module $\I_n$. 
The dynamics is given by modular multiplication by $\lambda$, 
and the study of periodicity involves asymptotic problems that 
are quite similar to those of the rational case
\cite{PercivalVivaldi,BartuccelliVivaldi,DegliEspostiIsola}.
The symbolic representation of hyperbolic toral automorphism 
was considered in \cite{KenyonVershik}, and is related to the 
codings corresponding to expansions to non-integral bases, 
the so-called {\it beta-expansions}
\cite{Sidorov,AkiyamaBorbeliBrunottePethoThuswaldner,%
AkiyamaBrunottePethoThuswaldner}.
We note that a module structure for periodic points is also 
found for certain dynamical systems over finite fields.
The set of $n$-cycles of a linearized polynomial\footnote{A 
polynomial over a finite field is linearized if all exponents 
are powers of the characteristic of the field.} 
is a vector space, that is, a module over a field 
\cite[lemma 4.2]{BatraMorton} 
(see also \cite{CohenHachenberger}).

Of a different nature are the codings related to renormalization,
which first appeared in dynamics in the study of persistence
of planar rotational motion under perturbation \cite{ShenkerKadanoff,MacKay}.  
In this case the scaling process is a second-order 
recursion over a suitable space of maps, and the 
symbols associated to it are the coefficients of the 
continued fraction expansion of the rotational angle.
As a result, the renormalization dynamics is eventually
periodic when the rotational angle is a quadratic irrational.
(These dynamical systems should not be confused with
the so-called rotational expansions \cite{Sidorov}.)

In this paper we consider a specific one-dimensional dynamical
system  ---an interval exchange transformation proposed by 
Arnoux and Yoccoz--- which features exact scaling, meaning that
the induced map (Poincar\'e section) on one of the intervals 
is a scaled version of the original map.
This map is defined over a cubic number field ---see
equation (\ref{eq:f}).
An interesting feature of this model is the existence of a symbol 
space which can be used to describe both the dynamics of the map 
and the dynamics of scaling. The dynamics of the map has zero
topological entropy, and its symbolic representation could be described 
as a generalized odometer \cite{LowensteinPoggiaspallaVivaldi:06}; 
by contrast, the dynamics of scaling is a shift with positive entropy.

This type of coding originates from the so-called 
{\it recursive tiling property\/} of a dynamical system:
a scale-invariant set of tiles exists, which cover the entire 
space under iteration, in a hierarchical fashion.
The symbol space consists of pairs of integers, which label a tile 
together with its transit time measured from a reference position.
These tilings have been used extensively in the study of piecewise 
isometric systems in two dimensions
\cite{KouptsovLowensteinVivaldi,LowensteinKouptsovVivaldi,%
LowensteinPoggiaspallaVivaldi}. 
It turns out that these symbol sequences posses a rich arithmetical
structure, which, to our knowledge, is largely unexplored.
In this respect, the Arnoux-Yoccoz model is an ideal object of study;
it illustrates how basic constructs developed elsewhere generalize, 
while remaining sufficiently simple for the numerical exploration of 
its asymptotics.

This paper is devoted to the study of the scaling dynamics 
---a subshift of finite type--- and in particular of its periodic orbits. 
This mapping leads to an arithmetic coding where the digits are algebraic 
integers rather than integers (a phenomenon already studied in 
unrelated contexts \cite{Gilbert,KenyonVershik}) and where the 
eventually periodic sequences coincide with the elements of the 
underlying field.

The system under study is defined in section \ref{section:Preliminaries},
where some key results in \cite{LowensteinPoggiaspallaVivaldi:06} are 
reviewed.
In section \ref{section:PeriodicPoints} we show that the periodic 
points belong to the cartesian product of a module over a finite ring, 
which depends on the period, and a finite set ---the {\it core region},
which is period-independent.
These two components represent the `fractional and integer parts'
of the periodic points, respectively.
The finiteness of the core region derives from the Pisot property 
of the scaling factor.
In section \ref{section:OrderingByPeriod}, we order the periodic points
by increasing period, and consider asymptotic properties of their 
integer and fractional parts. 
We prove that the set of periodic points which share the fractional
part with some other periodic point is infinite.
Numerical evidence, however, strongly suggests that this set has zero
density, which is expressed as conjecture \ref{conjecture:Multiplicity}.
By embedding the fractional parts of the cycles into the unit cube,
we study the relation between fractional and integer parts from
a geometrical viewpoint.
We collect some evidence that the phase space is organized in a 
highly nontrivial manner, which suggests directions for future 
investigations.

In section \ref{section:OrderingByDenominator}, we order the
periodic points by increasing `denominator', and relate their 
period to the ideal factorization of the denominator 
in the ring of integers of the underlying number field 
(proposition \ref{proposition:Order}).
Our experiments here investigate the asymptotic distribution 
of periods for large denominator, obtaining a picture
consistent with that of increasing period.

\vspace*{10pt}
\noindent
This research was supported by EPSRC grant No GR/S62802/01.

\section{Preliminaries}\label{section:Preliminaries}

Arnoux and Yoccoz \cite{ArnouxYoccoz}, in their study of 
pseudo-Anosov diffeomorphisms, introduced a family of 
interval-exchange transformations defined over algebraic 
number fields of arbitrarily large degree, with the property
of being renormalizable, namely of generating only finitely
many induced maps, up to scaling.
The simplest nontrivial member of this family corresponds
to a cubic field; the object of our study is the scale-invariant 
(exactly renormalizable) mapping $\rho$ induced by it on a sub-interval.
This mapping $\rho:[0,1)\to[0,1)$ is defined as follows 
\cite{ArnouxYoccoz,LowensteinPoggiaspallaVivaldi:06}
\begin{equation}\label{eq:Rho}
\rho(x)=x+\tau_j,\quad x\in\Omega_j
\hskip 40pt
\Omega_j=[\delta_j,\delta_{j+1})\qquad j=0,\ldots,N-1
\end{equation}
where $N=7$, the discontinuity points $\delta_N=1$ and translations 
$\delta_j$, $\tau_j$ are listed in table \ref{table:AYdata}, 
and $\lambda$ is the real root 
of the irreducible cubic polynomial
\begin{equation}\label{eq:f}
f(x)=x^3+x^2+x-1.
\end{equation}
\begin{table}[h]
$$
\begin{array}{c|c|c|c|l}
j&\delta_j&\tau_j&\nu_j&\quad p(j,t)\\
\hline
0&0&\lambda+\lambda^2&4&(0,6,3,6)\\
1&1-\lambda-\lambda^2&-1+3\lambda&13&(0,6,3,6,1,6,2,5,6,1,6,3,6)\\
2&1-2\lambda+\lambda^2&\lambda-\lambda^2&12&(0,6,3,6,1,6,2,5,6,2,4,6)\\
3&(3-4\lambda-\lambda^2)/2&-1+2\lambda+\lambda^2&8&(0,6,3,6,1,6,3,6)\\
4&(1-\lambda^2)/2&1-\lambda-\lambda^2&8&(0,6,4,5,6,2,4,6)\\
5&(-1+2\lambda+3\lambda^2)/2&\lambda-\lambda^2&12&(0,6,4,5,6,2,5,6,1,6,3,6)\\
6&\lambda&-\lambda&4&(0,6,4,6)
\end{array}
$$
\caption{\label{table:AYdata}
The data defining the cubic scale-invariant Arnoux-Yoccoz map $\rho$.
}
\end{table}
The map $\rho$ corresponds to the permutation 
$(0,1,2,3,4,5,6)\mapsto (6,2,4,3,5,1,0)$ of the $N$ disjoint 
intervals $\Omega_j$, which constitute a partition of $[0,1)$.
The resulting dynamical system has zero topological entropy, being
a piecewise isometry \cite{Buzzi}; in addition, it is uniquely 
ergodic, and hence, in particular, it has no periodic orbits.

The map $\rho$ is defined over the cubic number field 
\begin{equation}\label{eq:QLambda}
\Q(\lambda)=\{r_0+r_1\lambda+r_2\lambda^2\,|\,r_i\in\Q\}.
\end{equation}
If we restrict the coefficients $r_i$ in (\ref{eq:QLambda}) to 
integer values, we obtain the ring $\Z[\lambda]$ of all 
algebraic integers in $\Q(\lambda)$;
its group of units (invertible elements) 
has rank one, and is generated by $\lambda$ \cite[p.~519]{Cohen}. 
The unit $\omega=\lambda^3$ is important: in fact, scaling by $\omega$
conjugates the map $\rho$ to the induced map on the leftmost 
interval $\Omega_0=[0,1-\lambda-\lambda^2)=[0,\omega)$,
which can be verified by direct calculation.

To the above data we associate the sequence space $\V$
constituted by the sequences
$$
\sigma=\left((j_1,t_1), (j_2,t_2),\ldots\right)
\hskip 40pt
0\leq j_k <N, \quad 0\leq t_k < \nu_{k}
$$
where the $j_k$ are subject to the constraints
\begin{equation}\label{eq:PathConstraint}
j_k=p(j_{k+1},t_{k+1}),\quad k\geq 1
\end{equation}
and $p$ is the {\it path function\/} defined in table \ref{table:AYdata}.
Each sequence $\sigma\in\V$ satisfying (\ref{eq:PathConstraint}) 
defines a real number $x$ via the equations
\begin{equation}\label{eq:Point}
x(\sigma)=\sum_{i=1}^\infty d_i\omega^{i-1}
\hskip 40pt
d_i=\sum_{t=0}^{t_i-1}\tau_{p(j_i,t)}
\end{equation}
where the `digit' $d_i=d_i(j_i,t_i)$ belongs to a {\it finite\/} 
subset of the $\Z$-module $\M$ generated by the translations $\tau_k$. 
It is easy to see that $\M=\Z[\lambda]$; furthermore,
there are $\sum \nu_j=61$ valid symbols, corresponding to 25 
distinct values of the digits.

By construction, the point $x(\sigma)$ belongs to the unit 
interval. The set of points sharing a finite digit sequence 
$(d_1,\ldots,d_n)$ is a half-open interval, called a 
{\it tile\/} of level $n$: it is the intersection of 
translated and scaled versions of the original intervals $\Omega_j$.
For each integer $n$, the tiles of level $n$ constitute a 
partition of the unit interval: this is the {\it recursive tiling
property\/} of the mapping, which implies the surjectivity
of the function $\sigma\mapsto x(\sigma)$.
The representation (\ref{eq:Point}) may be regarded as a variant 
of the expansion of real numbers in non-integral bases ---the
so called $\beta$-expansions. Our basis $\omega$ is the reciprocal 
of a Pisot number\footnote{A Pisot number is an algebraic number greater 
than one, whose algebraic conjugates other than itself lie inside 
the complex unit circle.}, a property which is of importance here
(it underpins theorems 
\ref{theorem:EventuallyPeriodicCodes}, \ref{theorem:Decomposition}
below) as it is in the theory of $\beta$-expansions \cite{Akiyama,Sidorov}.
We emphasize again that the digits of our $\omega$-expansion are algebraic integers.


The map $\sigma\mapsto x(\sigma)$ is not injective, but the lack
of injectivity is easily controlled.  Whenever it fails, there 
are precisely two codes corresponding to the same point, and 
these codes are eventually 
periodic and entirely characterized by their periodic part.
This degeneracy can be removed by stipulating that the nested set of 
tiles specified by the code $\sigma$ must also contain its limit point $x$. 
If only such codes are allowed, the correspondence between points and codes 
is bi-unique. This amounts to excluding the following six periodic tails
\begin{equation}\label{eq:BadCodes}
(1,9)^\infty\quad (2,6)^\infty\quad (3,2)^\infty
\quad (4,6)^\infty\quad (5,3)^\infty\quad (6,1)^\infty.
\end{equation}
(The situation is analogous to the decimal representation of 
rational numbers: the tail $(9)^\infty$ is excluded to achieve
a unique representation.)

With the above stipulation on valid codes, the map
\begin{equation}\label{eq:Shift}
\gamma: [0,1)\to[0,1) \hskip 40pt x\mapsto (x-d_1)\omega^{-1}
\end{equation}
is conjugate to the left shift on $\V$.
This is a subshift of finite type, which we call
the {\it scaling dynamics\/} of the Arnoux-Yoccoz map.
The transitivity of the incidence matrix corresponding to the 
admissibility conditions (\ref{eq:PathConstraint}) can be verified
directly. It then follows \cite[p.~51]{KatokHasselblatt}, that 
$\gamma$ is mixing, and that the periodic orbits are dense.

The following result was proved in 
\cite{LowensteinPoggiaspallaVivaldi:06}.
\begin{theorem} \label{theorem:EventuallyPeriodicCodes}
A point $x$ of the unit interval belongs to
$\Q(\lambda)$ if and only if the code $\sigma(x)$ is eventually periodic.
\end{theorem}
This result calls to mind the familiar property of expansions in
integral bases, where the set of points characterized 
by the eventual periodicity of the digits is the field $\Q$.
To place this theorem into context, we note the general results of
Bertrand \cite{Bertrand} and Schmidt \cite{Schmidt} on eventually 
periodic greedy expansions in a Pisot base.

To study dynamics over $\Q(\lambda)$ we consider the decomposition
$$
\Q(\lambda)= \Xi + \Z[\lambda]
$$
where
$$
\Xi=\{\xi_0+\xi_1\lambda+\xi_2\lambda^2\,|\,\xi_i\in \Q \cap [0,1)\}.
$$
This prescription gives a unique representation of $x\in\Q(\lambda)$ as 
$x=\xi+\beta$, with $\xi\in\Xi$ and $\beta\in\Z[\lambda]$.
Accordingly, we define the projections to the two components
\begin{equation}\label{eq:Projections}
\Pi_f:\Q(\lambda)\to\Xi\qquad x\mapsto \xi,
\hskip 50pt
\Pi_i:\Q(\lambda)\to\B\qquad x\mapsto \beta.
\end{equation}
These components should be viewed, respectively, as the 
`fractional and integer parts' of the field elements. 
Note that $\Xi$ is a $\Z[\lambda]$-module, isomorphic to 
$\Q(\lambda)/\Z[\lambda]$; this fact will be of importance later.

In the rest of this paper we shall be concerned with the set of periodic points 
of $\gamma$, a $\gamma$-invariant set which belongs to $\Q(\lambda)\cap[0,1)$.

\section{Periodic points of scaling dynamics}\label{section:PeriodicPoints}

Consider a point $x$ with strictly periodic code 
$$
\sigma(x)=((j_1,t_1),\ldots,(j_n,t_n))^\infty.
$$
From theorem \ref{theorem:EventuallyPeriodicCodes} we have that 
$x\in\Q(\lambda)\cap [0,1)$, and a straightforward calculation shows that
\begin{equation}\label{eq:PeriodicPoint}
x=\frac{1}{1-\omega^n} \,\sum_{i=1}^n d_i\,\omega^{i-1}
\end{equation}
where the digits $d_i$ were defined in (\ref{eq:Point}).
The following result ---proved in
\cite{LowensteinPoggiaspallaVivaldi:06}---
is crucial to our investigation
\begin{theorem} \label{theorem:Decomposition}
There exists a finite subset $\B$ of $\Z[\lambda]$ such that
any $\gamma$-periodic point $x$ can be represented as
\begin{equation}\label{eq:Decomposition}
x\,=\,\xi+\beta 
\hskip 40pt \xi\in\Xi,\, \beta\in \B.
\end{equation}
Conversely, for every $\xi\in\Xi$, there exists $\beta\in\B$  such
that $x=\xi+\beta$ is periodic under $\gamma$.
\end{theorem}
The set $\B\in\Z[\lambda]$ is called the {\it core region.} 
The theorem says that the points which are periodic under the scaling 
map $\gamma$ can have only finitely many distinct integer parts. 
Exact computations give the estimates
\begin{equation}\label{eq:Bounds}
31\leq\#\B\leq 348.
\end{equation}
The upper bound was established using a technique 
developed in \cite{LowensteinPoggiaspallaVivaldi:06}. 
(In that reference, an upper bound of 225 was obtained for 
the projection of the core region to a two-dimensional lattice.) 
The lower bound will be established below 
---see table \ref{table:CoreRegion}.

Theorem \ref{theorem:Decomposition} outlines the basic structure
of the periodic points of $\gamma$, yet important questions 
concerning their integer and fractional parts are not addressed.
Are there infinitely many fractional parts $\xi\in\Xi$ shared 
by more than one periodic point?
Is it true that for all $\beta$ in the core region
there are infinitely many periodic points with integer 
part equal to some $\beta$?
Furthermore, the above questions can be refined to formulate
statements on the density of the relevant infinite sets.

To investigate these issues, for $\beta\in\B$ we consider the set
$\Xi(\beta)$ of all $\xi\in\Xi$ for which $\xi+\beta$ is periodic.
Using the projections defined in (\ref{eq:Projections}), we see that  
\begin{equation}\label{eq:XiBeta}
\Xi{(\beta)}=\Pi_f\circ\Pi^{-1}_\beta(\beta)
\end{equation}
and from theorem \ref{theorem:Decomposition} we know that
\begin{equation}\label{eq:XiPartition}
\bigcup_{\beta\in\B}{\Xi}{(\beta)}=\Xi.
\end{equation}
Next, for $\xi\in\Xi$ we define the {\it multiplicity\/} $\kappa(\xi)$ 
to be the number of $\gamma$-periodic points with fractional part 
equal to $\xi$. From theorem \ref{theorem:Decomposition} we know 
that $1\leq \kappa(\xi)\leq \#\B$, and we are interested in the 
set $\partial\Xi$ of those $\xi$ that correspond to more than 
one periodic point. This set is given by
\begin{equation}\label{eq:BoundarySet}
\partial \Xi=\left\{\xi\in\Xi\,:\,\kappa(\xi)>1\right\}\,=\,
\bigcup_{\beta,\beta'\in\B\atop\beta\neq\beta'}\left(\Xi(\beta)\cap\Xi(\beta')\right).
\end{equation}
For reason that will become clear later, we call $\partial \Xi$
the {\it boundary set.}

We begin our study of periodic points by seeking an analogue 
of equation (\ref{eq:Lattice}).
Let $M_n$ be the smallest positive integer which is
divisible by $1-\omega^n$ in $\Z[\lambda]$, that is,
\begin{equation}\label{eq:M}
(1-\omega^n)\alpha_n=M_n
\end{equation}
for $M_n$ minimal and some $\alpha_n\in\Z[\lambda]$.
From equation (\ref{eq:PeriodicPoint}) one sees that $M_nx\in\Z[\lambda]$,
that is,
$$
\mbox{Fix}(\gamma^n)\subset \frac{1}{M_n}\Z[\lambda].
$$
Equivalently, $M_n$ is the least common multiple of the denominators 
of all $n$-periodic points. 
In practise, $M_n$ is computed efficiently with the extended Euclid's 
algorithm.

By construction, $M_n$ is a divisor of the norm $N(1-\omega^n)$ of 
$1-\omega^n$ (the norm of an algebraic number is the product of
its algebraic conjugates).
While the value of $N(1-\omega^n)$, is well-behaved ---see section 
\ref{section:OrderingByPeriod}, equations (\ref{eq:Norm},\ref{eq:NoIn})---
the value of $M_n$ features large fluctuations.
In this respect, we note the factorization
\begin{equation}\label{eq:Cyclotomic}
1-\omega^n=-\prod_{d|n}C_{3d}(\lambda)
 \prod_{d|n\atop 3\not \hskip 1.5pt |\,d}C_{d}(\lambda)
\end{equation}
where $C_d(x)$ is the $d$-th cyclotomic polynomial \cite[p 37]{Niven}.
This formula is established starting from the polynomial identity 
$x^n-1=\prod_{d|n}C_d(x)$, recalling that $\omega=\lambda^3$, and then
relating the divisors of $3n$ to those of $n$.
The formula provides only a partial factorization of $1-\omega^n$,
because the irreducibility of the cyclotomic polynomial $C_d(x)$ 
in $\Z[x]$ does not imply that $C_d(\lambda)$ are primes in 
$\Z[\lambda]$ (which is a principal ideal domain \cite[p.~519]{Cohen}). 

For example, the identities
$$
\omega-1=C_1(\lambda)C_3(\lambda)=(\lambda-1)\lambda^{-1}
\hskip 40pt
(1-\lambda)^3\lambda^{-5}=2 
$$
together with the fact that $\lambda$ is a unit, show that the 
denominator of a periodic point always has the irreducible 
factor $1-\lambda$, which is a prime divisor of 2.
From the above observation, it follows that $M_n$ is even.
Likewise, $C_6(\lambda)$ is a divisor of 7, and so if $n$
is even, $M_n$ is divisible by 14, etc. 

The periodic orbits of the scaling dynamics can be ordered in two
ways, namely by increasing period $n$, or by increasing denominator $m$
(the choice of sub-orderings for given $n$ or $m$ are not important here). 
As mentioned above, these orderings lead to different perspectives
(see also~\cite[section 6]{ChirikovVivaldi}), and they will be
considered in sections \ref{section:OrderingByPeriod} and 
\ref{section:OrderingByDenominator}, respectively.
In analogy with the rational case, whenever $1-\omega^n$ has a 
large number of divisors, cancellation becomes likely in equation 
(\ref{eq:PeriodicPoint}), leading to $n$-cycles with denominator 
much smaller than $M_n$. As the denominator is increased, these 
are the cycles that are found first.

Each ordering gives asymptotic information about the sets 
$\Xi(\beta)$ and $\partial\Xi$ defined above. 
This will be our main object of study.
 
\section{Ordering by period}\label{section:OrderingByPeriod}
Let the period $n$ be fixed.  With reference to equations
(\ref{eq:Projections}) and (\ref{eq:Decomposition}), we define
\begin{equation}\label{eq:DecompositionII}
\I_n=\Pi_f\left(\frac{1}{1-\omega^n}\Z[\lambda]\right)
\hskip 30pt
\I_n^\prime=\Pi_f\left(\mbox{Fix}(\gamma^n)\right)
\hskip 30pt
\B_n=\Pi_i\left(\mbox{Fix}(\gamma^n)\right).
\end{equation}
From equation (\ref{eq:PeriodicPoint}) it follows that the
fractional part of the $n$-cycles is contained in $\I_n$,
and hence $\I_n^\prime\subset\I_n$. By construction,
the fractional parts in the (possibly empty) residual 
set $\I_n\setminus\I_n'$ belong to cycles whose period 
is a multiple of $n$.

To determine the cardinality of $\I_n$, we must refine the 
decomposition (\ref{eq:Decomposition}). 
Let $m$ be a positive integer, and let $\Xi_{m}$ be the set 
of $\xi$-points with denominator $m$ namely
\begin{equation}\label{eq:Xim}
\Xi_m=\Pi_f\left(\frac{1}{m}\Z[\lambda] \right)
\end{equation}
which is again a $\Z[\lambda]$-module, and is isomorphic 
(via the map $\varphi: x\mapsto m x$) to the module
$\Z[\lambda]/m\Z[\lambda]$, which has $m^3$ elements.

Consider now the algebraic integer $\alpha_n$ defined in equation
(\ref{eq:M}). 
The map $\varphi$ sends $\I_n$ to the ideal generated by 
$\alpha_n$ in $\Z[\lambda]/M_n\Z[\lambda]$ (now regarded as a ring), 
which yields the equation 
$$
\#\Xi_{M_n}=\#\I_n\times |N(\alpha_n)|.
$$
Using the multiplicativity of the norm, and the fact that 
$\#\Xi_M=N(M)=M^3$, we obtain the formula
\begin{equation}\label{eq:Norm}
\#\I_n=|N(1-\omega^n)|
\end{equation}
which, together with equation (\ref{eq:M}), shows that
$\#\I_n$ is a multiple of $M_n$.
Explicit computation gives
\begin{equation}\label{eq:NoIn}
\#\I_n
 =\omega^{-n}-\omega^{n}-2\left(\sqrt{\omega^{-n}}-\sqrt{\omega^{n}}\right)\cos(n\theta)
\end{equation}
where
$$
\cos(\theta)=\frac{1}{2}\sqrt{\omega}\,(5-\omega).
$$

To compute $\mbox{\# Fix}(\gamma^n)$ we consider the 
{incidence matrix}
$$
A=\left(\begin{array}{ccccccc}
1&0&0&1&0&0&2\\
1&2&1&2&0&1&6\\
1&1&2&1&1&1&5\\
1&1&0&2&0&0&4\\
1&0&1&0&2&1&3\\
1&1&1&1    &1&2&5\\
1&0&0&0&1&0&2\end{array}\right)
$$
where $A_{i,j}$ is the number of times that the first-return orbit of the 
scaled interval $\omega\Omega_i$ visits the interval $\Omega_j$, 
obtained from the path function data of table \ref{table:AYdata}.
The characteristic polynomial of $A$ factors into irreducibles as 
$(x-1)(x^3-5x^2+7x-1)(x^3-7x^2+5x-1)$, whose real roots are 
$1$, $\omega=\lambda^3$, and $\omega^{-1}$, respectively.
Accounting for the six forbidden period-1 codes listed in 
(\ref{eq:BadCodes}),
we find
\begin{equation}\label{eq:NoFixn}
\mbox{\# Fix}(\gamma^n)=\mbox{Tr}A^n-6
 =\omega^{-n}+\omega^{n}+2(\sqrt{\omega^{-n}}+\sqrt{\omega^{n}})\cos(n\theta)-5.
\end{equation}
Expressions (\ref{eq:NoIn}) and (\ref{eq:NoFixn}) have the same 
leading term $\omega^{-n}$ (recall that $|\omega|<1$), and
hence $\mbox{\# Fix}(\gamma^n)\sim\mbox{\#}\I_n$. Furthermore
\begin{equation}\label{eq:Difference}
\mbox{\# Fix}(\gamma^n)-\mbox{\#}\I_n=4\sqrt{\omega^{-n}}\cos(n\theta)+2\omega^n-5
\,\sim\, 4\sqrt{\omega^{-n}}\cos(n\theta).
\end{equation}
This formula shows that for some values of $n$, $\# \mbox{Fix}(\gamma^n)$ 
is greater than $\# \I_n$, which implies that there exist distinct periodic 
points with the same fractional part $\xi$.
In fact, the asymptotic expression (\ref{eq:Difference}) 
shows that this must happen infinitely often, and that the
population of points with multiplicity greater than 1 
is maximal when the period is a denominator of the continued 
fraction expansions of $\theta/2\pi$.
Recalling the definition (\ref{eq:BoundarySet}) of the boundary
set, we have the following result.
\begin{proposition} \label{proposition:BoundarySet}
The boundary set $\partial\Xi$ is infinite.
\end{proposition}
Using induction, it is possible to prove that the
periodic codes
$$
((1,4)^k,(3,4),(1,2))^\infty
\hskip 40pt
((2,9)^k,(4,5),(2,10))^\infty
\qquad
k=0,1,\ldots
$$
(for $k=0$ the symbol is deleted) correspond to points of 
period $k+2$ on the boundary set. Now, the fixed points
$(1,4)^\infty$ and $(2,9)^\infty$ also belong to $\partial\Xi$,
as easily verified.
This shows that the boundary set contains points of any period,
which is stronger than proposition \ref{proposition:BoundarySet}.
We shall not produce this proof here.

When $\# \mbox{Fix}(\gamma^n)<\# \I_n$, the $\xi$-values do not 
exhaust the whole of $\I_n$, and formula (\ref{eq:Difference})
again shows that this must happen for infinitely often.
Since $\#{\I'}_n\#\B \geq \#\mbox{Fix}(\gamma^n)$, we obtain 
\begin{equation}\label{eq:Bound}
\frac{\#\I_n'}{\#\I_n}\geq 
\frac{\#\mbox{Fix}(\gamma_n)}{\#\I_n \#\B}
\,\to\, \frac{1}{\#\B}\qquad\mbox{as}\quad  n\to\infty.
\end{equation}
This bound is far from optimal. To see this,
we define
$$
{\I_n''}=\{\xi\in\I_n'\,:\,\kappa(\xi)>1\}.
$$
Constructing explicitly all $n$-periodic points for $n\leq 14$, 
gives the data displayed in table \ref{table:CyclesData}.

\begin{table}[h]
\hskip 100pt\vtop{\baselineskip=12pt\halign{
\quad \hfil $#$&\qquad \hfil$#$&\qquad \hfil$#$\hfil\qquad&\hfil$#$\quad\cr
n & \#{\I''}_n & \#{\I''}_n/\mbox{\# Fix}(\gamma^n) &\# \B_n\cr
\noalign{\vskip 3pt\hrule\vskip 5pt}
1&2   &0.2857& 7\cr
2&14  &0.3111&17\cr
3&38  &0.1467&21\cr
4&86  &0.05595&22\cr
5&182 &0.01939&24\cr
6&374 &0.006440&24\cr
7&758 &0.002100& 27\cr
8&1534&0.0006831& 30\cr
9&3170&0.0002268&30\cr
10&-&-&30\cr
11&-&-&30\cr
12&-&-&30\cr
13&-&-&30\cr
14&-&-&31\cr
}}
\caption{\label{table:CyclesData} Multiplicity and core region data
from orbits of period $n$.}
\end{table}

The multiplicity data are limited to the range $n\leq 9$, 
because the comparison of all fractional parts requires 
that these quantities be stored.
The most significant finding, to be used in conjunction with the
asymptotic formulae (\ref{eq:NoIn}) and (\ref{eq:NoFixn}), is
that distinct periodic points typically have distinct fractional 
parts; namely, the set of cycles with multiplicity 1 has full density. 
Furthermore, the ratio $\#{\I''}_n/\mbox{\# Fix}(\gamma^n)$ 
appears to be tending rather rapidly to zero, with a scaling 
ratio close to $1/3$.
In the computations to be described in section 
\ref{section:OrderingByDenominator} we observe the same 
phenomenon when the cycles are ordered by increasing 
denominator, leading to the following conjecture.
\begin{conjecture}\label{conjecture:Multiplicity}
The following holds
$$
\lim_{n\to\infty}\frac{\#\{\xi\in {\I'}_n \,:\,\kappa(\xi)=1 \}}
{\# \I_n} \,=\,
\lim_{m\to\infty}\frac{\#\{\xi\in\Xi_m\,:\,\kappa(\xi)=1 \}}{m^3} \,=\,1.
$$
\end{conjecture}
The validity of this conjecture would imply that, as $n\to\infty$
$$
\#\I_n^\prime\,\sim\, \#\I_n\,\sim\,\omega^{-n}
$$
to be compared with the bound (\ref{eq:Bound}).

Regarding the core region data of table \ref{table:CyclesData},
the lower bound of 31 is substantially smaller than 
the upper bound of 348 given in (\ref{eq:Bounds}).
With reference to equation (\ref{eq:Decomposition}), we display 
in table \ref{table:CoreRegion} some elements 
$\beta=m_0+m_1\lambda+m_2\lambda^2$ of the core region, 
represented as integer triples $(m_0,m_1,m_2)$.
Each value of $\beta$ has a probability $\mu(\beta)$, defined as 
the limiting density of the $\xi$-values that correspond to it.
Approximate values of these probabilities may be computed in two
different ways, corresponding to ordering by period (the second
and third columns in the table) and by denominator (the fourth column).
In the former case, the data was computed from the {14}-cycles 
(130399019341 data points), which project to the 31 points
of the core region listed in the first column, respectively. 
The data for period 13 are also shown, to give an idea of convergence.
The probability $\mu'$ computed using all cycles with `denominator'
$m\leq 200$ is displayed in the third column (373112717 data 
points, which project to a subset of 28 core region points)
---see the next section. 
While it is conceivable that the core region could contain more 
points than those displayed in table \ref{table:CoreRegion}, 
the stability of the above figures suggests that they
are a reliable approximation to the densities.

\begin{table}
\hskip 60pt\vtop{\baselineskip=12pt\halign{
\quad $#$&\hfil $#$&\hfil $#$&\hfil $#$&\qquad$#$\hfil&\quad $#$\hfil& \quad $#$\quad\hfil\cr
\noalign{\hskip 40pt$\beta$\hskip 50pt $\mu_{13}$\hskip 60pt $\mu_{14} $\hskip 62pt $\mu'$
          \vskip 3pt\hrule\vskip 5pt}
(& 0,& 0,& 0)& 0.21275 & 0.21275 & 0.21566\cr
(& 0,&-1,&-1)& 0.18921& 0.18921 & 0.18811\cr
(& 0,& 0,&-1) & 0.13777& 0.13777 & 0.13688\cr
(& 1,&-2,&-1) & 0.095625 & 0.095625 & 0.096146\cr
(&-1,& 0,& 0) & 0.077156 & 0.077155 & 0.075485\cr
(&-1,& 0,& 1) & 0.062212 & 0.062213 & 0.062157\cr
(& 0,&-2,& 0) & 0.57953 & 0.057952 & 0.057331\cr
(& 0,&-1,& 1) & 0.029309& 0.029308 & 0.029370\cr
(& 1,&-3,& 0) & 0.026513 & 0.026514 & 0.026326\cr
(& 1,&-2,& 0) & 0.026416& 0.026416 & 0.027329\cr
(& 0,&-2,&-1) & 0.022015& 0.022014 & 0.020936\cr
(& 0,&-2,& 1) & 0.013304 & 0.013304 & 0.013345\cr
(& 0,&-1,& 0) & 0.011008& 0.011009 & 0.010908\cr
(& 0,& 0,& 1) & 0.0088052 & 0.0088053 & 0.0093162\cr
(&-1,& 1,& 0) & 0.0079820 & 0.0079830 & 0.0083504\cr
(&-1,& 1,& 1) & 0.0068145& 0.0068150 & 0.0069532\cr
(& 0,& 1,& 0) & 0.0056184& 0.0056188 & 0.0058062\cr
(& 1,&-1,&-1) & 0.0052609& 0.0052609 & 0.0056897\cr
(& 1,&-3,&-1) & 0.0020413& 0.0020418 & 0.0019191\cr
(&-1,& 0,&-1) & 0.0011704& 0.0011707 & 0.0010408\cr
(&-1,&-1,& 1) & 0.00075012& 0.00075025 & 0.00065573\cr
(& 1,&-4,& 0) & 0.00020058& 0.00020069 & 0.00018132\cr
(& 0,&-3,& 0) & 7.5098\cdot10^{-5} & 7.5092\cdot10^{-5} & 5.2630 
\cdot10^{-5}\cr
(& 1,&-2,& 1) & 2.1778\cdot10^{-5} & 2.1757\cdot10^{-5} & 3.2518 
\cdot10^{-5}\cr
(& 1,&-3,& 1) &3.7165\cdot10^{-6}& 3.7439\cdot10^{-6} & 7.8207\cdot10^ 
{-6}\cr
(& 0,&-1,&-2) & 3.6249\cdot10^{-6}& 3.6509\cdot10^{-6} & 1.1444 
\cdot10^{-6}\cr
(& 1,& 0,&-1) &3.3849 \cdot10^{-6}& 3.4058\cdot10^{-6} & 7.8207 
\cdot10^{-6}\cr
(& 1,& -1,& -2) &2.6184 \cdot 10^{-7}&2.6498\cdot 10^{-7} & 0\cr
(& 0,& -3,& 1) & 2.6184\cdot 10^{-7}&2.6498\cdot 10^{-7} & 4.3151 
\cdot 10^{-7}\cr
(& -1,& 1,& -1) & 9.1552\cdot 10^{-8}&9.3050\cdot 10^{-8} & 0\cr
(&-1,&-1,& 0) & 0 & 1.2873 \cdot10^{-11} & 0\cr
}}
\caption{\label{table:CoreRegion} The elements $\beta$ of the core
region $\B$, with three distinct estimates $\mu_{13}$, $\mu_{14}$ and  
$\mu'$ of the associated densities. The first two were computed using  
all cycles up to periods 13 and 14, respectively, and the third using  
all cycles with denominator not exceeding 200.} The floating-point  
numbers represent exact rationals rounded off to 5-digit precision.
\end{table}

So far we have expressed probabilistic information in terms 
of densities; we close this section by discussing related
questions of measure. To this end, we consider the embedding 
of the fractional parts into the unit cube
\begin{equation}\label{eq:Embedding}
\varphi:\Xi\to [0,1)^3\hskip 40pt
r_0+r_1\lambda+r_2\lambda^2\mapsto (r_0,r_1,r_2),\qquad 0\leq r_i<1.
\end{equation}
This construct depends on the particular basis chosen for $\Q(\lambda)$.
However, a change of basis corresponds to a unimodular transformation,
which, if we identify the unit cube with the 3-torus, is continuous
and volume preserving. So any property of the (closure of) embedded 
objects which is topological or concerns three-dimensional Lebesque 
measure will be independent of the basis.

It is natural to consider sequences of points that are 
uniformly distributed in the cube; for instance, the 
sequence of fractional parts of periodic points with 
increasing period, or with increasing denominator.
Consider now the embedding $\varphi(\partial\Xi)$ of the 
boundary set $\partial\Xi$, defined in equation (\ref{eq:BoundarySet}).
Recalling that $\#\mbox{Fix}(\gamma^n)\sim\omega^{-1}=\lambda^{-3}$
(see equation (\ref{eq:NoFixn})), the scaling law observed 
in table \ref{table:CyclesData} is consistent with 
the closure of $\varphi(\partial\Xi)$ having zero three-dimensional 
Lebesque measure.
Now, the points of the boundary set are precisely the points 
with multiplicity greater than one; so any uniformly distributed 
sequence will have the property that the density of the elements 
with multiplicity greater than 1 is zero. 
Thus the vanishing of the measure of the boundary set in
the cube implies conjecture \ref{conjecture:Multiplicity}. 

We have examined the arrangements of the boundary set in the cube,
which gives some support to this conjecture, although the 
available data are not quite conclusive.
In fact, there is an even stronger property that should be
considered, namely that the core region localizes in 
the unit cube; this means that the measure of the closure 
of the sets $\varphi(\Xi(\beta))$ add up to unity, and 
represent the probabilities $\mu(\beta)$ estimated in 
table \ref{table:CoreRegion}.
(Loocalization techniques have proved very useful in the 
study of certain lattice maps \cite{LowensteinVivaldi,BosioVivaldi}.)

These problems deserve further investigation.

\section{Ordering by denominator}\label{section:OrderingByDenominator}
Let $m$ be a positive integer. We consider the periodic
orbits of the map $\gamma$ having denominator $m$.
According to theorem \ref{theorem:Decomposition}
and equation (\ref{eq:Xim}), these points belong to the 
set $\Xi_m +\B$. We are interested in the determination of 
their period.

In section \ref{section:OrderingByPeriod} we have seen that 
scaling by $m$ maps $\Xi_m$ into the ring 
$\Z[\lambda]/m\Z[\lambda]$, which for our purpose is a 
more convenient representation.
If $\alpha, \beta\in\Z[\lambda]$, we denote by $(\alpha)$, 
$(\alpha,\beta)$, etc., the ideals they generate. 
Thus $(1)=\Z[\lambda]$.

Let $\xi$ be a point in $\Z[\lambda]/(m)$. 
The {\it order\/} $t(\xi)$ of $\xi$ is defined as 
\begin{equation}\label{eq:Order}
t(\xi)=\min\{k\geq 1:\omega^k\xi\equiv \xi\,\mod{(m)}\}.
\end{equation}
From theorem \ref{theorem:Decomposition}, we have that every 
$\xi\in\Xi$ is the fractional part of at least one periodic 
point; it follows that the order of $\xi$ must divide the 
corresponding period. 
In fact, from conjecture \ref{conjecture:Multiplicity}
it would follow that the period of a periodic point is
typically equal to the order of its fractional part.
Now $\# \Xi_m=m^3$, which gives the crude bound $t(\xi)\leq m^3$. 

To obtain a sharper bound for $t$, we must consider the ideal 
factorization of $m$ in $\Z[\lambda]$. We describe it in terms of 
the factorization of $f(x)$ modulo a prime $p$.
Since $\Z[\lambda]$ is the ring of all algebraic 
integers in $\Q(\lambda)$, the discriminant of $f(x)$,
which is equal to $-44$, does not have any spurious 
prime divisor, so that the factorization of $f(x)$ modulo 
$p$ describes the ideal factorization of $p$, without 
exceptions \cite[theorem 27]{Marcus}.
The primes that divide the discriminant are $p=2,11$; they 
{\it ramify,} that is, $f(x)$ has multiple roots modulo $p$.
For all the other primes $p$, the polynomial $f(x)\mod{p}$ has 
distinct factors. A prime $p$ is {\it inert,} {\it splits,} 
or {\it splits completely\/}, respectively, if $f(x)$ decomposes 
into the product of 1,2, or 3 irreducible factors, respectively.
Each possibility occurs infinitely often, and indeed with probability 
1/3, 1/2, or 1/6, respectively.
This follows from Cebotarev's density theorem and the fact that the 
Galois group of $f(x)$ is $S_3$ \cite[p.~129]{PohstZassenhaus}.

The following result relates the order of a point $\xi$ to
the ideal factorization of its denominator. 
\begin{proposition}\label{proposition:Order}
Let $\xi\in\Z[\lambda]$, and let $m$ be a positive integer.
Then $t(\xi)$ is a divisor of $T$, where $T=T(m)$ is computed as follows:
\begin{enumerate}
\item If  $gcd(n,m)=1$, then $T(nm)={\mbox lcm}(T(n),T(m))$.
\item If $m=p^e$ with $p$ prime and $e\geq 1$, we have 

\hfil\vtop{\hsize=10cm
\halign{\quad #\hfil &\qquad $#$ \quad\hfil\cr
$\quad p$&\quad T(p^e)\cr
\noalign{\vskip 3pt \hrule\vskip 5pt}
ramif. & \cases{p^{e+1}(p-1) & if $p=2$\cr
                p^e(p-1) & if $p=11$\cr
}\cr
\noalign{\vskip 10pt}
inert & 
\cases{p^{e-1}(p^2+p+1)   & if $p\equiv 0,2\mod 3$\cr
       p^{e-1}(p^2+p+1)/3 & if $p\equiv 1\mod 3$\cr
}\cr
\noalign{\vskip 10pt}
splits & \, p^{e-1}(p^2-1)/3\cr
\noalign{\vskip 10pt}
splits c.& 
   \cases{p^{e-1}(p-1) & if $p\equiv 2\mod 3$\cr
          p^{e-1}(p-1)/3 & if $p\equiv 1\mod 3.$\cr
}\cr
}}\hfil
\end{enumerate}
\end{proposition}
\proof
With reference to the the congruence (\ref{eq:Order}), we see that
if $\xi$ and $m$ are coprime (meaning that the ideal identity 
$(\xi,m)=(1)$ holds), then upon division by $\xi$, we obtain
$\omega^t\equiv 1\mod{(m)}$, independent of $\xi$, which
shows that $t(\xi)$ is the multiplicative order of $\omega$ in the 
ring $\Z[\lambda]/(m)$. Thus $t(\xi)$ is a divisor of the 
order of the multiplicative group of this ring. From the decomposition
$$
\frac{\Z[\lambda]}{(m)}\simeq 
 \bigoplus_i \frac{\Z[\lambda]}{J_i}
$$
where $J_i$ are the pairwise coprime ideals divisors of $(m)$, we have, in
particular, that $T(m)$ can be computed as the least common multiple
of its value at each summand.

If $\xi$ and $m$ are not coprime, then division by $\xi$
leads to a congruence modulo the ideal $(m)/(\xi,m)$. 
The multiplicative group of the corresponding finite ring is a divisor 
of that considered above, so the same estimates apply. 

Next we consider primary factors $m=p^e$, beginning with the
case $e=1$.
We denote by $\F_{p^k}$ the field with $p^k$ elements,
and by $\F_{p^k}^*$ its multiplicative group.
We shall use the fact that $\omega=\lambda^3$ is 
a unit of norm 1, and that the prime $p=3$ is inert.

If $p$ is inert, we have the ring isomorphism
$$
\frac{Z[\lambda]}{(p)}\simeq \F_{p^3}. 
$$
The Galois group of $\F_{p^3}$ over $\F_p$ is cyclic of order 3, and is 
generated by the Frobenius automorphism ${\cal F}:\,\alpha\mapsto\alpha^p$. 
The orbits of ${\cal F}$ are the sets of algebraic conjugates.
Denoting again by $\omega$ its reduction to $\F_{p^3}$, we find
$$
\omega^{1+p+p^2}={\cal F}^0(\omega)\,{\cal F}^1(\omega)\,{\cal F}^2(\omega)=1.
$$
The last equality derives from the fact that the terms in 
the product are the algebraic conjugates of $\omega$, and 
since they are distinct, their product is the norm of $\omega$ modulo $p$.
So the order of $\omega$ divides $p^2+p+1$. If $p\equiv 1\mod{3}$,
then $p^2+p+1\equiv 0\mod{3}$, and $\omega$ has a cube root in $\F_{p^3}$, 
which is also a unit of norm 1. In this case
the order of $\omega$ divides $(p^2+p+1)/3$.

If $p$ splits then
$$
\frac{Z[\lambda]}{(p)}\simeq \F_{p}\oplus \F_{p^2}.
$$ 
The order of $\F_p^*$ divides
that of $\F_{p^2}^*$, which is equal to $p^2-1$.
Because $p\not=3$, we have $p^2-1\equiv 0\mod{3}$,
hence $\omega$ always has a cube root in $\F_{p^2}$.

If $p$ splits completely, then 
$$
\frac{Z[\lambda]}{(p)}\simeq \F_p\oplus\F_p\oplus\F_p.
$$ 
The order of the three reductions of $\omega$ to each finite 
field is a divisor of the order of $\F_p^*$, which is $p-1$. 
If $p\equiv 2\mod{3}$, there are no further constraints on $T(p)$. 
If $p\equiv 1\mod{3}$, then $\omega$ has a cube root in $\F_p$, 
and hence the subgroup it generates has index at least 3.

It remains to deal with the ramified primes 2 and 11. 
The prime 2 is totally ramified
$$
f(x)\equiv(x+a)^3\mod{p}\hskip 40pt p=2,\, a=1
$$
leading to the ideal factorization $(p)={P}^3=(p,\lambda+a)^3$. 
The ideal ${P}$ has $p^2$ incongruent points modulo 
$(p)$, and hence the multiplicative group of $\Z[\lambda]/(p)$
has order $p^3-p^2=p^2(p-1)$.

The prime $11$ is partially ramified. We find
$$
f(x)\equiv(x+a)^2(x+b)\mod{p}\hskip 40pt p=11,\, a=3,\, b=6
$$
which corresponds to the ideal factorization, 
$(p)={P}_1^2{P}_2=(p,\lambda+a)^2(p,\lambda+b)$.
Hence
$$
\frac{Z[\lambda]}{(p)}\simeq R \oplus \F_{p}.
$$
Where $R$ is a ring with $p^2$ elements, $p$ of which are 
not invertible. Thus its multiplicative group has order
$p^2-p=p(p-1)$, while that of $\F_p$ has order $p-1$,
hence the result.

Finally, for $e>1$, we apply a standard lifting argument.
Let $k$ be the largest integer for which $T(p^k)=T(p)$. 
(Such an integer $k$ exists and is effectively computable.)
That is,
$$
\omega^{T(p)}=1+p^k\beta
$$
where $\beta\in\Z[\lambda]$ is such that $(\beta, p)=(1)$.
Then, using the binomial theorem, we have
$$
\omega^{T(p)p}=1+p^{k+1}\beta'
$$
where again $(\beta',p)=(1)$. An easy induction on $e$ shows 
that the order increases regularly by a factor $p$ at each step.
Thus 
$$
T(p^e)=T(p)\,p^{\max(0,e-k)}.
$$

This completes the proof.
\endproof

In the totally ramified case, the $\xi$-dynamics is 
quite regular. From the fact that $(2)=2_1=(2,\lambda+1)^3$, and 
$\omega=1-\lambda-\lambda^2$, we have 
$\omega-1=-\lambda(\lambda+1)\in{2_1}$.
Thus, for any $\xi\in\Z[\lambda]$, we have that
$\omega\xi \equiv \xi \mod{2_1}$, and hence
every orbit of $\omega$ modulo $(2^e)$ consists of congruent 
points modulo ${2_1}$. So, if the period is maximal, orbits
are cosets of a lattice.

The theorem states that $t(\xi)$ is a divisor of $T$, but
gives no information as to the actual value of $T/t$.
The exact value of $t(\xi)$ is (essentially) determined by 
the order of the image of $\omega$ in various finite fields,
which cannot be computed in non-polynomial time. One would 
expect that the occurrence of a given ratio $T/t$ will be 
accompanied by large fluctuations and by a limiting probability, 
which should be maximal for $T/t=1$. As a rule, these are very 
hard problems, of the kind mentioned in the introduction
in connection with Artin's conjecture.

Proposition \ref{proposition:Order} generalizes, with appropriate 
modifications, to any cubic field. What is essential is that 
$\omega$ be a unit: the fact that our $\omega$ is the third 
power of a fundamental unit has merely brought about some 
specialization in the formulae, expressed via congruences. 
The presesence of the split case in proposition 
\ref{proposition:Order} is a consequence of the fact 
that the Galois group of $f(x)$ is the symmetric group.
If the Galois group of $f(x)$ is cyclic, the split case
does not occur, while the other formulae in the proposition 
remain the same. In this case the Cebotarev's densities 
for the allowed factorizations become 1/3 for primes 
splitting completely, and 2/3 for inert primes.

We have constructed all periodic points with denominator 
not exceeding 200. They have the form $x=\xi+\beta$
with $\beta\in\B$ and $\xi\in\cup_{m\leq 200} \Xi_m$,
where $\Xi_m$ was defined in (\ref{eq:Xim});
we obtain 373112717 points in total.
Now, every fractional part $\xi$ is represented in this set,
and the total number of fractional parts is given by
$$
\#\bigcup_{m\leq 200} \Xi_m=
 \sum_{m\leq 200}\bigg(\sum_{d|m}d^3\mu(m/d)\bigg)=373111960
$$
where $\mu$ is the M\"obius function. 
(The divisor sum gives the number of lattice points in $\Xi_m$ 
that do not belong to any lattice of smallest index.)
The difference between the number of periodic points 
and the number of their fractional parts gives the size 
of the boundary set restricted to these denominators
(to match the figures, one has to take into account the fact 
that two $\xi$-values correspond to three and four $\beta$-values, 
respectively).
In our data set the boundary points are one in $10^6$.
These findings agree with the corresponding data for 
orbits with increasing period, and underpin our 
conjecture \ref{conjecture:Multiplicity}.

\end{document}